      \documentclass[11pt]{amsart}
      \usepackage[mathscr]{eucal}
      \usepackage{amsmath,amsfonts}
\def\rg{\hbox to 30pt{\rightarrowfill}}
\def\lg{\hbox to 30pt{\leftarrowfill}}

      \parskip\smallskipamount
          \newtheorem{theorem}{Theorem}[section]
      
      \newtheorem{proposition}[theorem]{Proposition}
      \newtheorem{corollary}[theorem]{Corollary}

      \newtheorem{remark}[theorem]{Remark}
      \makeatletter
      \@addtoreset{equation}{section}
      \makeatother

\newcommand{\HH}{{\mathbb H}}
      \newcommand{\BB}{{\mathbb B}}
      \newcommand{\CC}{{\mathbb C}}
\newcommand{\KK}{{\mathbb K}}

      \newcommand{\DD}{{\mathbb D}}
      
      \newcommand{\FF}{{\mathbb F}}
      \newcommand{\TT}{{\mathbb T}}

      \newcommand{\cA}{{\mathcal A}}

      \newcommand{\cD}{{\mathcal D}}

      \newcommand{\cH}{{\mathcal H}}
      \newcommand{\cK}{{\mathcal K}}
      
      \newcommand{\cM}{{\mathcal M}}
      \newcommand{\cN}{{\mathcal N}}
      \newcommand{\cP}{{\mathcal P}}

      \newcommand{\cW}{{\mathcal W}}

      \newdimen\expt
      \def\boxit#1{\setbox0\hbox{$\displaystyle{#1}$}
            \hbox{\lower.4\expt
       \hbox{\lower3\expt\hbox{\lower\dp0
            \hbox{\vbox{\hrule height.4\expt
       \hbox{\vrule width.4\expt\hskip3\expt
            \vbox{\vskip3\expt\box0\vskip2\expt}%
       \hskip3\expt\vrule width.4\expt}\hrule height.4\expt}}}}}}
      
\begin{document}
       \pagestyle{myheadings}
      \markboth{ Gelu Popescu}{ Operator theory  on noncommutative varieties II }

      \title [  Operator theory  on noncommutative varieties II   ]
      {  Operator theory  on noncommutative varieties II }
        \author{Gelu Popescu}
      \date{September 7, 2005}
      \thanks{Research supported in part by an NSF grant}

      \address{Department of Mathematics, The University of Texas
      at San Antonio \\ San Antonio, TX 78249, USA}
      \email{\tt gelu.popescu@utsa.edu}

      \begin{abstract}
An $n$-tuple of operators $T:=[T_1,\ldots, T_n]$ on a Hilbert space
$\cH$ is called a $J$-constrained row contraction if  \ $
T_1T_1^*+\cdots + T_nT_n^*\leq I_\cH$    and
$$
f(T_1,\ldots, T_n)=0,\quad f\in J,
$$
where $J$ is a WOT-closed two-sided ideal of the noncommutative analytic Toeplitz
 algebra $F_n^\infty$ and $f(T_1,\ldots, T_n)$ is defined using the
  $F_n^\infty$--functional calculus for row contractions.
  We show that
the {\it constrained characteristic function} $\Theta_{J,T}$
associated with $J$ and $T$ is  a  complete unitary invariant for
$J$-constrained  completely non-coisometric (c.n.c.) row
contractions.  We  also provide a model   for  this class of  row
contractions in terms of the constrained characteristic functions.
In particular,   we obtain  a  model theory   for
 $q$-commuting c.n.c. row contractions.

 \end{abstract}
      \maketitle
      \section*{Introduction}

In \cite{Po-constrained}, we developed a dilation theory
on noncommutative  varieties determined by  row contractions $T:=[T_1,\ldots, T_n]$
subject to constraints such as
$$
p(T_1,\ldots, T_n)=0,\quad p\in \cP,
$$
where $\cP$ is a set of noncommutative polynomials. In this setting, the model $n$-tuple is the universal row contraction $[B_1,\ldots, B_n]$ satisfying the same constraints as $T$, which turns out to be   the   maximal  constrained piece of the $n$-tuple $[S_1,\ldots, S_n]$ of left creation operators on the full Fock space on $n$ generators.
We obtained a Beurling type theorem characterizing the invariant subspaces under each operator $B_1\otimes I_\cH,\ldots, B_n\otimes I_\cH$, and Wold type decompositions for $*$-representations of the $C^*$-algebra $C^*(B_1,\ldots, B_n)$
generated  by $B_1,\ldots, B_n$ and the identity.  The constrained dilation and model theory is based on a class of
{\it constrained Poisson kernels} associated with $T$  and representations of the Toeplitz algebra $C^*(B_1,\ldots, B_n)$.

Following the classical Sz.-Nagy--Foia\c s model theory for a single contraction \cite{SzF-book} and  the  multivariable noncommutative dilation theory \cite{Bu}, \cite{Fr}, \cite{Po-models}, \cite{Po-isometric}, \cite{Po-charact}, we introduced in \cite{Po-constrained}
a {\it constrained characteristic function} $\Theta_{J,T}$  associated with any constrained row contraction $T$. It turned out that,
  for  constrained {\it pure} row contractions,
the constrained  characteristic function  is  a  complete unitary invariant.   We also showed that the curvature invariant and Euler characteristic
asssociated with a Hilbert module generated by  an arbitrary
(resp.~commuting) row contraction $T$ can be expressed only in terms of the
(resp.~constrained) characteristic function of $T$.
This paper is a continuation of \cite{Po-constrained}.
We further investigate the constrained characteristic function in several variables.

      Let $H_n$ be an $n$-dimensional complex  Hilbert space with orthonormal
      basis
      $e_1$, $e_2$, $\dots,e_n$, where $n\in\{1,2,\dots\}$.        We consider the full Fock space  of $H_n$ defined by
      $$F^2(H_n):=\bigoplus_{k\geq 0} H_n^{\otimes k},$$
      where $H_n^{\otimes 0}:=\CC 1$ and $H_n^{\otimes k}$ is the (Hilbert)
      tensor product of $k$ copies of $H_n$.
      Define the left creation
      operators $S_i:F^2(H_n)\to F^2(H_n), \  i=1,\dots, n$,  by
      $$
       S_i\varphi:=e_i\otimes\varphi, \quad  \varphi\in F^2(H_n).
      $$
      The    noncommutative analytic Toeplitz algebra   $F_n^\infty$
        and  its norm closed version,
        the noncommutative disc
       algebra  $\cA_n$,  were introduced by the author   in \cite{Po-von} (see also  \cite{Po-funct}).        $F_n^\infty$  is the algebra of left multipliers of
      the Fock space $F^2(H_n)$  and  can be identified with
       the
        weakly closed  (or $w^*$-closed) algebra generated by the left creation
      operators
         $S_1,\dots, S_n$    and the identity.
           The noncommutative disc algebra $\cA_n$ is
          the  norm closed algebra generated by
          the same operators.
           When $n=1$, $F_1^\infty$
         can be identified
         with $H^\infty(\DD)$, the algebra of bounded analytic functions
          on the open unit disc. The  noncommutative analytic Toeplitz algebra $F_n^\infty$ can be viewed as a
           multivariable noncommutative
          analogue of $H^\infty(\DD)$.

Let $T:=[T_1,\ldots, T_n]$  be a {\it completely non-coisometric}
(c.n.c.) row contraction (see Section 1 for notation) with  $T_i\in
B(\cH)$, the algebra of all bounded linear operators
        operators on a Hilbert space $\cH$. Given a
WOT-closed two-sided ideal $J$ of the noncommutative analytic
Toeplitz algebra $F_n^\infty$, we say that $T$ is   a
$J$-constrained row contraction if
$$
f(T_1,\ldots, T_n)=0,\quad f\in J,
$$
where  $f(T_1,\ldots, T_n)$ is defined using the
$F_n^\infty$-functional calculus for c.n.c. row contractions
\cite{Po-funct}. In Section 1, we  present some results concerning
constrained Poisson transforms associated with $J$-constrained row
contractions. More about   noncommutative  Poisson kernels and
Poisson transforms on $C^*$-algebras generated by isometries can be
found in \cite{Po-poisson}, \cite{APo}, \cite{Po-tensor},
\cite{Po-curvature},  \cite{Po-unitary}, and \cite{Po-constrained}.

 The {\it constrained characteristic function} associated with  a $J$-constrained row contraction $T$ was introduced in
 \cite{Po-constrained} as   a   multi-analytic operator
\, $ \Theta_{J,T}:\cN_J\otimes \cD_{T^*}\to \cN_J\otimes \cD_T $\,
uniquely defined by the formal Fourier representation
$$ -I_{\cN_J}\otimes T+
\left(I_{\cN_J}\otimes \Delta_T\right)\left(I_{{\cN_J}\otimes \cH}-\sum_{i=1}^n W_i\otimes T_i^*\right)^{-1}\\
\left[W_1\otimes I_\cH,\ldots, W_n\otimes I_\cH
\right] \left(I_{\cN_J}\otimes \Delta_{T^*}\right),
$$
where $W_1,\ldots, W_n$ are the constrained right creation operators associated with $J$
(see Section 1 for notation).

In Section 2,   we show that the constrained characteristic function  is  a complete unitary
invariant
  for  the class of   constrained c.n.c. row contractions.
We  also provide a model   for  this class of  row contractions in terms of the constrained characteristic functions.
  All the results of this paper  apply, in particular,
to c.n.c. row contractions subject to constraints such as $
p(T_1,\ldots, T_n)=0$,\ $ p\in\cP, $ where $\cP$ is a set of
noncommutative polynomials.

In  particular,   we obtain a model theory for commuting c.n.c. row
contractions. The characteristic function of a commuting row
contraction $T:=[T_1,\ldots, T_n]$ is the compression to the
symmetric Fock space of the noncommutative characteristic function
introduced in \cite{Po-charact}. As shown in \cite{Po-constrained}
(using \cite{Arv}), it can be identified with the operator-valued
analytic function on the open unit ball  of $\CC^n$, given by
$$
\Theta_{J_c,T}(z):= -T+\Delta_T(I-z_1T_1^*-\cdots -z_nT_n^*)^{-1}
[z_1I_\cH,\ldots, z_nI_\cH]\Delta_{T^*}, \quad z=(z_1,\ldots,
z_n)\in \BB_n,
$$
In this particular setting, the characteristic function was proved
to be a complete unitary invariant for pure row contractions in
\cite{BES} and, independently, by the author
 in \cite{Po-constrained}. We should mention that, in the
 commutative case, there is a model theory (see e.g. \cite{AE} and
 references therein) for operator  tuples $T$ satisfying positivity
 conditions of type $\frac{1}{K}(T,T^*)\geq 0$, where $K$ is a
 reproducing kernel associated with certain domains in $\CC^n$. It
 will be interesting to see if there is a ``constrained-version'' of
 all this work.

After the completion of this paper, we received a preprint from
T.~Bhattacharyya, J.~Eschmeier, and J.~Sarkar \cite{BES2}, and  also
noticed a very recent paper of C.~Benhida, and D.~Timotin \cite{BT}.
Both papers  deal with the characteristic function of a commuting
row contraction and there is some overlap, in this particular
setting, with Section 2 of our paper.
  However, our results concerning the model theory for constrained row contractions
   are more general and the proofs are based on noncommutative dilation theory, rather
    then
   reproducing  kernel Hilbert space techniques.

\bigskip

\section{Constrained Poisson kernels}\label{constr.poisson}

 In this section, we present   some results concerning
the constrained Poisson kernels associated with  completely non-coisometric  row
 contractions. These results are needed in   Section 2.

We need to recall from
      \cite{Po-multi}, \cite{Po-von},  \cite{Po-funct},  \cite{Po-analytic}, and \cite{Po-central}
       a few facts
       concerning multi-analytic   operators on Fock
      spaces.
We denote by $\FF_n^+$  the unital free semigroup on $n$ generators
      $g_1,\dots,g_n$, and the identity $g_0$.
      The length of $\alpha\in\FF_n^+$ is defined by
      $|\alpha|:=k$ if $\alpha=g_{i_1}g_{i_2}\cdots g_{i_k}$, and
      $|\alpha|:=0$ if $\alpha=g_0$.
       If $T_1,\dots,T_n\in B(\cH)$,   define       %
      $T_\alpha :=  T_{i_1}T_{i_2}\cdots T_{i_k}$
      if $\alpha=g_{i_1}g_{i_2}\cdots g_{i_k}$, and
      $T_{g_0}:=I_\cH$.  Similarly, we denote $e_\alpha:=
e_{i_1}\otimes\cdots \otimes  e_{i_k}$ and $e_{g_0}:=1$.
         We say that
       a bounded linear
        operator
      $M$ acting from $F^2(H_n)\otimes \cK$ to $ F^2(H_n)\otimes \cK'$ is
       multi-analytic
      if
      \begin{equation*}
      M(S_i\otimes I_\cK)= (S_i\otimes I_{\cK'}) M\quad
      \text{\rm for any }\ i=1,\dots, n.
      \end{equation*}
       We can associate with $M$ a unique formal Fourier expansion
       $    M(R_1,\ldots, R_n):= \sum_{\alpha \in \FF_n^+}
       R_\alpha \otimes \theta_{(\alpha)}$,
       for some operators $\theta_{(\alpha)}\in B(\cK,\cK')$,
      where $R_i:=U^* S_i U$, \ $i=1,\ldots, n$, are the right creation
      operators
      on $F^2(H_n)$ and
      $U$ is the (flipping) unitary operator on $F^2(H_n)$ mapping
        $e_{i_1}\otimes e_{i_2}\otimes\cdots\otimes e_{i_k}$ into
       $e_{i_k}\otimes\cdots\otimes e_{i_2}\otimes e_{i_1} $.
       Since  the operator $M$ acts like its Fourier representation on ``polynomials'',
  we will identify them for simplicity.
      The set of  all multi-analytic operators in
      $B(F^2(H_n)\otimes \cK,
      F^2(H_n)\otimes \cK')$  coincides  with
      $R_n^\infty\bar \otimes B(\cK,\cK')$,
      the WOT-closed operator space generated by the spatial tensor product, where
      $R_n^\infty=U^* F_n^\infty U$.
A multi-analytic operator is called inner if it is an isometry, and outer if it has dense range.

Now let $J\neq F_n^\infty$ be  a WOT-closed two-sided ideal of  the
noncommutative analytic Toeplitz algebra $F_n^\infty$. Define the
subspaces of  the full Fock space $F^2(H_n)$ by setting
$$
\cM_J:=\overline{JF^2(H_n)}\quad \text{and}\quad \cN_J:=F^2(H_n)\ominus \cM_J.
$$
   Define the {\it constrained  left} (resp.~{\it right}) {\it creation operators}
   by setting
$$B_i:=P_{\cN_J} S_i|_{\cN_J} \quad \text{and}\quad W_i:=P_{\cN_J} R_i|_{\cN_J},
\quad i=1,\ldots, n.
$$
Let $\cW(B_1,\ldots, B_n)$ be the $w^*$-closed algebra generated by
$B_1,\ldots, B_n$ and the identity. We proved in \cite{APo}
that
$$
\cW(B_1,\ldots, B_n)=P_{\cN_J}F_n^\infty |_{\cN_J}=\{f(B_1,\ldots, B_n):\
f(S_1,\ldots, S_n)\in F_n^\infty\},
$$
where,  according to the $F_n^\infty$-functional calculus for c.n.c.
row contractions \cite{Po-funct},
$$f(B_1,\ldots, B_n)=\text{\rm SOT-}\lim\limits_{r\to 1}f(rB_1,\ldots, rB_n).
$$
Note that if $\varphi\in J$, then $\varphi(B_1,\ldots, B_n)=0$.
 An operator $M\in B(\cN_J\otimes \cK,\cN_J\otimes \cK')$ is called multi-analytic with respect to
the constrained shifts $B_1,\ldots, B_n$ if
$$
M(B_i\otimes I_{\cK})=(B_i\otimes I_{\cK'})M,\quad i=1,\ldots, n.
$$
If in addition $M$ is partially isometric, then we call it inner. If $M$ has dense range, it is called outer.
We recall from \cite{Po-central} that the set of all multi-analytic  operators with respect to
 $B_1,\ldots, B_n$ coincides  with
$$
\cW(W_1,\ldots, W_n)\bar\otimes  B(\cK,\cK')=P_{\cN_J\otimes \cK'}[R_n^\infty\bar\otimes B(\cK,\cK')]|_{\cN_J\otimes \cK},
$$
and a similar result holds for  the algebra $\cW(B_1,\ldots, B_n)$.

\smallskip

Now, let us review  (see \cite{Po-poisson}) some basic properties for  noncommutative
Poisson transforms associated with row contractions $T:=[T_1,\ldots, T_n]$,
\ $T_i\in B(\cH)$. For  each $0<r\leq 1$, define the defect operator
$\Delta_{T,r}:=(I-r^2T_1T_1^*-\cdots -r^2 T_nT_n^*)^{1/2}$.
The Poisson  kernel associated with $T$ is the family of operators
$$
K_{T,r} :\cH\to F^2(H_n)\otimes \overline{\Delta_{T,r}\cH}, \quad  0<r\leq 1,
$$
defined by

\begin{equation}\label{Poiss}
K_{T,r}h:= \sum_{k=0}^\infty \sum_{|\alpha|=k} e_\alpha\otimes r^{|\alpha|}
\Delta_{T,r} T_\alpha^*h,\quad h\in \cH.
\end{equation}
When $r=1$, we denote $\Delta_T:=\Delta_{T,1}$ and $K_T:=K_{T,1}$.
The operators $K_{T,r}$ are isometries if $0<r<1$, and
\begin{equation}
\label{K*K}
 K_T^*K_T=I_\cH- \text{\rm SOT-}\lim_{k\to\infty}
\sum_{|\alpha|=k} T_\alpha T_\alpha^*.
\end{equation}
This shows that $K_T$ is an isometry if and only if $T$ is a {\it pure} row
 contraction (\cite{Po-isometric}),
i.e.,
\begin{equation}\label{KK}
\text{\rm SOT-}\lim_{k\to\infty} \sum_{|\alpha|=k} T_\alpha T_\alpha^*=0.
\end{equation}
A key property of the Poisson kernel
is that
\begin{equation}\label{eq-ker}
K_{T,r}(r^{|\alpha|} T_\alpha^*)=(S_\alpha^*\otimes I)K_{T,r}\qquad
\text{ for all } 0<r\leq 1,\ \alpha\in \FF_n^+.
\end{equation}

 When $T$ is a completely non-coisometric (c.n.c.) row-contraction, i.e.,
 there is no $h\in \cH$, $h\neq 0$, such that
 \begin{equation} \label{cnc}
 \sum_{|\alpha|=k}\|T_\alpha^* h\|^2=\|h\|^2
 \quad \text{\rm for all } \ k=1,2,\ldots,
 \end{equation}
an  $F_n^\infty$-functional calculus was developed  in \cite{Po-funct}.
 We showed that if $f=\sum\limits_{\alpha\in \FF_n^+} a_\alpha S_\alpha$ is
 in $F_n^\infty$, then
 \begin{equation}\label{fc}
 \Gamma_T(f)=f(T_1,\ldots, T_n):=
 \text{\rm SOT-}\lim_{r\to 1}\sum_{k=0}^\infty
  \sum_{|\alpha|=k} r^{|\alpha|} a_\alpha T_\alpha
\end{equation}
exists and $\Gamma_T:F_n^\infty\to B(\cH)$ is a WOT-continuous
completely contractive homomorphism.

Let $J\neq F_n^\infty$ be a WOT-closed two-sided ideal of
$F_n^\infty$, and let $T:=[T_1,\ldots, T_n]$, \ $T_i\in B(\cH)$, be
a row contraction. The {\it constrained Poisson  kernel} associated
with $J$ and $T$ is the operator $K_{J,T}:\cH\to \cN_J\otimes
\overline{\Delta_T \cH}$  defined by $ K_{J,T}:=(P_{\cN_J}\otimes
I_{ \overline{\Delta_T \cH}}) K_T, $ where $K_T$ is the Poisson
kernel defined by relation \eqref{Poiss} (case $r=1$).

The following theorem provides extensions of   some results from
\cite{Po-unitary} for   constrained pure row contractions.

\begin{theorem}\label{KJT}
Let $J\neq F_n^\infty$ be a WOT-closed two-sided  ideal of
$F_n^\infty$ and let $T:=[T_1,\ldots, T_n]$, $T_i\in B(\cH)$, be a
c.n.c. row contraction such that $$ \varphi(T_1,\ldots, T_n)=0,
\qquad \varphi\in J. $$ Then
$$
K_{J,T} f(T_1,\ldots, T_n)^*=(f(B_1,\ldots, B_n)^*\otimes I_{\overline{\Delta_T\cH}})K_{J,T}
$$
for any $f(B_1,\ldots, B_n)\in \cW(B_1,\ldots, B_n)$, where $K_{J,T}$ is the constrained Poisson kernel
associated with $J$ and $T$.
\end{theorem}
\begin{proof}

 Due to  relation \eqref{eq-ker} (case
$r=1$)
  we have
 \begin{equation}\label{ker-pol}
 K_T^*(p(S_1,\ldots, S_n)\otimes I_\cH) =p(T_1,\ldots, T_n)K_T^*
 \end{equation}
 for any   polynomial $p(S_1,\ldots, S_n)$. According to \cite{Po-funct},
 if $f(S_1,\ldots, S_n):= \sum\limits_{k=0}^\infty
 \sum\limits_{|\alpha|=k} a_\alpha S_\alpha$
 is
 in $F_n^\infty$, then, for any $0<r<1$,
 $f_r(S_1,\ldots, S_n):= \sum\limits_{k=0}^\infty \sum\limits_{|\alpha|=k}
 r^{|\alpha|}a_\alpha S_\alpha$  is in the noncommutative disc algebra $\cA_n$.
 Since
  $$
 \lim_{m\to \infty} \sum_{k=0}^m \sum_{|\alpha|=k}
  r^{|\alpha|} a_\alpha S_\alpha = f_r(S_1,\ldots, S_n)
 $$
  in the norm topology, relation
   \eqref{ker-pol} implies
  \begin{equation*}
 K_T^*(f_r(S_1,\ldots, S_n)\otimes I_\cH) =f_r(T_1,\ldots, T_n)K_T^*
 \end{equation*}
for any $f(S_1,\ldots, S_n)\in F_n^\infty$ and $0<r<1$.
Since $T$ is a c.n.c. row contraction and $S:=[S_1,\ldots, S_n]$ is a pure row contraction, we can
 use the $F_n^\infty$-functional calculus.
  We recall that the map $A\mapsto A\otimes I$ is SOT-continuous on bounded sets of $B(F^2(H_n))$ and, due to the noncommutative von Neumann inequality \cite{Po-von} (see \cite{vN} for the classical case),  we have $\|f_r(S_1,\ldots, S_n)\|\leq \|f(S_1,\ldots, S_n)\|$. Therefore,  we can
 take  $r\to 1$ in
 the above equality and
  obtain
 \begin{equation}\label{ker-pol3}
K_T^*(f(S_1,\ldots, S_n)\otimes I_{\overline{\Delta_T\cH}}) =f(T_1,\ldots, T_n)K_T^*
 \end{equation}
 for any $f(S_1,\ldots, S_n)\in F_n^\infty$, where $f(T_1,\ldots, T_n)$
 is defined
 by formula \eqref{fc}.

 Now, according to   relation \eqref{ker-pol3}, we have
\begin{equation}
 \label{ker-p}
\left<(\varphi(S_1,\ldots, S_n)^*\otimes I)K_Th,1\otimes k\right>=
 \left<K_T\varphi(T_1,\ldots, T_n)^*h,1\otimes k\right>
 \end{equation}
 for any $\varphi(S_1,\ldots, S_n)\in F_n^\infty$,\  $h\in \cH$, and $k\in {\overline{\Delta_T\cH}}$.
 Note that if $\varphi(S_1,\ldots, S_n)\in J$, then
 $\varphi(T_1,\ldots, T_n)=0$, and relation \eqref{ker-p} implies
 $\left<K_Th, \varphi\otimes k\right>=0$ for any $h,k\in \cH$.
 Taking into account the definition of $\cM_J$, we deduce that
 \begin{equation}
 \label{KT-sub}
 K_T(\cH)\subseteq \cN_J\otimes \cH.
 \end{equation}
  This shows that the constrained Poisson kernel satisfies the relation
 \begin{equation}
 \label{KT}
 K_{J,T}h=\left(P_{ \cN_J}\otimes I_{\overline{\Delta_T \cH}}\right) K_Th=K_T h, \quad h\in \cH.
 \end{equation}

 Since $J$ is a left ideal of $F_n^\infty$, $\cN_J$ is an invariant subspace
 under each operator $S_1^*,\ldots, S_n^*$ and therefore
  $B_\alpha=P_{ \cN_J}S_\alpha|\cN_J$,
 $\alpha\in \FF_n^+$.
 Since $[B_1,\ldots, B_n]$ is a pure row contraction, we can use
the $F_n^\infty$-functional calculus to deduce that

 \begin{equation}
 \label{fnn}
 f(B_1,\ldots, B_n)=P_{ \cN_J}f(S_1,\ldots, S_n)|\cN_J
 \end{equation}
 for any $f(S_1,\ldots, S_n)\in F_n^\infty$.
Taking into account relations \eqref{ker-pol3},
 \eqref{KT}, and \eqref{fnn},
we obtain
\begin{equation*}
 \begin{split}
K_{J,T}f(T_1,\ldots, T_n)^*&=\left(P_{ \cN_J}\otimes I_{\overline{\Delta_T\cH}}\right)
[f(S_1,\ldots, S_n)^*\otimes I_{\overline{\Delta_T \cH}}]\left(P_{ \cN_J}\otimes I_{\overline{\Delta_T\cH}}\right)K_T\\
&=
 \left[\left(P_{ \cN_J}f(S_1,\ldots, S_n)|
 \cN_J\right)^*\otimes I_{\overline{\Delta_T \cH}}\right] K_{J,T}\\
&= \left[f(B_1,\ldots, B_n)^*\otimes I_{\overline{\Delta_T \cH}}\right]K_{J,T}.
 \end{split}
 \end{equation*}
Therefore, we have
\begin{equation}\label{KFFK}
 K_{J,T}f(T_1,\ldots, T_n)^*=\left[f(B_1,\ldots, B_n)^*\otimes I_{\overline{\Delta_T \cH}}\right]K_{J,T}
 \end{equation}
 for any $f(B_1,\ldots, B_n)\in \cW(B_1,\ldots, B_n)$.
This completes the proof.
\end{proof}

\begin{corollary}\label{Co}
Let $J\neq F_n^\infty$ be a WOT-closed two-sided  ideal of $F_n^\infty$ and let $T:=[T_1,\ldots, T_n]$, $T_i\in B(\cH)$, be a c.n.c. row contraction such that
$$
\varphi(T_1,\ldots, T_n)=0, \qquad \varphi\in J.
$$
Then
$$
K_{J,T} T_i^*=(B_i^*\otimes I_{\overline{\Delta_T\cH}})K_{J,T}, \quad i=1,\ldots, n,
$$
 and
\begin{equation}\label{J-pure}
K_{J,T}^* K_{J,T}=
I_\cH- \text{\rm SOT-}\lim_{k\to\infty} \sum_{|\alpha|=k}
T_\alpha T_\alpha^*,
\end{equation}
where $K_{J,T}$ is the constrained Poisson kernel associated with
$T$ and $J$.

Moreover, if a map $\Psi: \cW(B_1,\ldots, B_n)\to B(\cH)$ satisfies
the relation
\begin{equation}\label{intert}
\Psi(f) K_{J,T}^*=K_{J,T}^* (f\otimes I),\quad f\in \cW(B_1,\ldots,
B_n),
\end{equation}
then
$$\Psi(f)=f(T_1,\ldots, T_n),\quad f\in \cW(B_1,\ldots, B_n).
$$
\end{corollary}
 \begin{proof} The first part of the theorem follows easily from Theorem \ref{KJT}.
 Relation \eqref{J-pure} is a consequence of \eqref{K*K} and \eqref{KT-sub}.  Due to relation \eqref{J-pure}, if $T$ is a c.n.c. row contraction  then  $K_{J,T}$ is a one-to-one operator. Consequently, Theorem \ref{KJT} and relation \eqref{intert} imply
$\Psi(f)=f(T_1,\ldots, T_n)$, $ f\in \cW(B_1,\ldots, B_n)$. This
completes the proof.
\end{proof}
We remark that if $\cP$ is  a family of noncommutative polynomials  in $S_1,\ldots, S_n$ and  $T:=[T_1,\ldots, T_n]$ is  an arbitrary
row contraction such that
$p(T_1,\ldots, T_n)=0$, $p\in \cP$, then one can prove that
$T$ is a c.n.c. row contraction if and only if  $K_{J,T}$ is a one-to-one operator, where $J$ is the WOT-closed two-sided  ideal of $F_n^\infty$ generated  by $\cP$.

  \bigskip

      \section{Constrained characteristic functions} \label{Characteristic}

In this section, we show that the constrained characteristic function  is  a complete unitary
invariant
  for  the class of constrained c.n.c. row contractions.
We  also provide a model   for  this class of  row contractions in terms of the constrained characteristic functions.
  All the results of this section  apply, in particular,
to c.n.c. row contraction subject to constraints such as
$$
p(T_1,\ldots, T_n)=0,\quad p\in\cP,
$$
where $\cP$ is a set of noncommutative polynomials.

The characteristic  function associated with an arbitrary row
contraction $T:=[T_1,\ldots, T_n]$, \ $T_i\in B(\cH)$, was
introduced in \cite{Po-charact} (see \cite{SzF-book} for the
classical case $n=1$) and it was proved to be  a complete unitary
invariant for completely non-coisometric (c.n.c.) row contractions.
Using the characterization of multi-analytic operators on Fock
spaces (see \cite{Po-analytic}, \cite{Po-tensor}), one can easily
see that the characteristic  function  of $T$ is  a multi-analytic
operator
$$
\Theta_T:F^2(H_n)\otimes \cD_{T^*}\to F^2(H_n)\otimes \cD_T
$$
with the formal Fourier representation
\begin{equation*}
\begin{split}
  -I_{F^2(H_n)}\otimes T+
\left(I_{F^2(H_n)}\otimes \Delta_T\right)&\left(I_{F^2(H_n)\otimes \cH}-\sum_{i=1}^n R_i\otimes T_i^*\right)^{-1}\\
&\left[R_1\otimes I_\cH,\ldots, R_n\otimes I_\cH
\right] \left(I_{F^2(H_n)}\otimes \Delta_{T^*}\right),
\end{split}
\end{equation*}
where $R_1,\ldots, R_n$ are the right creation operators on the full Fock space $F^2(H_n)$.
 Here,  we need to clarify some notations since some of them are different from those considered in \cite{Po-charact}.
The defect operators  associated with a row contraction $T:=[T_1,\ldots, T_n]$
are
\begin{equation*}
\Delta_T:=\left( I_\cH-\sum_{i=1}^n T_iT_i^*\right)^{1/2}\in B(\cH) \quad \text{ and }\quad \Delta_{T^*}:=(I-T^*T)^{1/2}\in B(\cH^{(n)}),
\end{equation*}
while the defect spaces are $\cD_T:=\overline{\Delta_T\cH}$ and
$\cD_{T^*}:=\overline{\Delta_{T^*}\cH^{(n)}}$, where
$\cH^{(n)}$ denotes the direct sum of $n$ copies of $\cH$.
   We proved \cite{Po-constrained} that
\begin{equation}
\label{fa-ca}
I_{F^2(H_n)\otimes \cD_T}-\Theta_T \Theta_T^*=K_T K_T^*,
\end{equation}
where $K_T$ is the Poisson kernel associated with $T$.

Let $J\neq F_n^\infty$ be a WOT-closed two-sided ideal of the noncommutative analytic Toeplitz algebra $F_n^\infty$.  In \cite{Po-constrained}, we defined the {\it constrained characteristic  function} associated with a $J$-constrained c.n.c.  row contraction
$T:=[T_1,\ldots, T_n]$, \ $T_i\in B(\cH)$,
to be   the  multi-analytic operator (with respect to the constrained shifts $B_1,\ldots, B_n$)
$$
\Theta_{J,T}:\cN_J\otimes \cD_{T^*}\to \cN_J\otimes \cD_T
$$
defined by the formal Fourier representation
$$ -I_{\cN_J}\otimes T+
\left(I_{\cN_J}\otimes \Delta_T\right)\left(I_{{\cN_J}\otimes \cH}-\sum_{i=1}^n W_i\otimes T_i^*\right)^{-1}\\
\left[W_1\otimes I_\cH,\ldots, W_n\otimes I_\cH
\right] \left(I_{\cN_J}\otimes \Delta_{T^*}\right).
$$
Taking into account that $\cN_J$ is  a co-invariant subspace under
$R_1,\ldots, R_n$, we  have
\begin{equation}\label{rest}
\begin{split}
\Theta_{T}^*(\cN_J\otimes\cD_T)
&\subseteq \cN_J\otimes \cD_{T^*}\ \text{  and  }\\
P_{\cN_J\otimes \cD_{T}}\Theta_{T}|\cN_J\otimes \cD_{T^*}&=\Theta_{J,T}.
\end{split}
\end{equation}
Let us remark that   the above definition
of the constrained characteristic function  makes sense  when
$T:=[T_1,\ldots, T_n]$ is an arbitrary    $J$-constrained row contraction
 and $J$ is a WOT-closed two-sided ideal of $F_n^\infty$ generated by a
 family of polynomials (see \cite{Po-constrained}).

The next result was obtained in \cite{Po-constrained} for
WOT-closed two-sided ideal of
$F_n^\infty$ generated by polynomials. Here, we have an extension of that result.

\begin{theorem}\label{J-factor}
Let $J\neq F_n^\infty$ be a WOT-closed two-sided ideal of   $F_n^\infty$ and
let  $T:=[T_1,\ldots, T_n]$, \ $T_i\in B(\cH)$,  be a $J$-constrained c.n.c. row
 contraction.
Then
\begin{equation}\label{J-fa}
I_{\cN_J\otimes \cD_T}-\Theta_{J,T}\Theta_{J,T}^*=K_{J,T}K_{J,T}^*,
\end{equation}
where
$\Theta_{J,T}$ is  the  constrained characteristic function  of   $T$ and $K_{J,T}$ is  the corresponding
constrained Poisson kernel.
\end{theorem}

\begin{proof}

Due to relation \eqref{KT-sub}, we have  $\text{\rm range}\,
K_T\subseteq \cN_J\otimes \overline{\Delta_T\cH}$. Taking the
compression of relation \eqref{fa-ca} to the subspace $\cN_J\otimes
\cD_T\subset F^2(H_n)\otimes \cD_T$, we obtain
$$
I_{\cN_J\otimes \cD_T}-P_{\cN_J\otimes \cD_T}\Theta_T
\Theta_T^*|\cN_J\otimes \cD_T=
P_{\cN_J\otimes \cD_T}K_TK_T^*|\cN_J\otimes \cD_T.
$$
Using relation
 \eqref{rest}  and that $W_i^*=R_i^*|\cN_J$, \ $i=1,\ldots, n$,  we deduce \eqref{J-fa}. The proof is complete.
\end{proof}

Now, we present a model for    constrained c.n.c. row contractions  in
terms of the constrained  characteristic functions.
We recall that two row contractions  $T:=[T_1,\ldots, T_n]$, $T_i\in B(\cH)$, and  $T'=:[T_1',\ldots, T_n']$, $T_i'\in B(\cH')$, are equivalent if there is a unitary operator $U:\cH\to \cH'$ such that $UT_i=T_i'U$ for any $i=1,\ldots, n$.

\begin{theorem}\label{model}
Let $J\neq F_n^\infty$ be a WOT-closed two-sided ideal of
$F_n^\infty$ and $T:=[T_1,\ldots, T_n]$ be a  c.n.c.  row contraction such that
$$
\varphi(T_1,\ldots, T_n)=0,\quad \varphi\in J.
$$
Then $T:=[T_1,\ldots, T_n]$ is unitarily equivalent to the constrained row contraction $\TT:=[\TT_1,\ldots, \TT_n]$ on
the Hilbert space
\begin{equation*}
\HH_{J,T}:=\left[\left(\cN_J\otimes \cD_T\right)\oplus \overline{\Delta_{J,T}(\cN_J\otimes \cD_{T^*})}\right]
\ominus\left\{\Theta_{J,T}f\oplus \Delta_{J,T}f:\ f\in \cN_J\otimes \cD_{T^*}\right\},
\end{equation*}
where $\Delta_{J,T}:= \left(I-\Theta_{J,T}^* \Theta_{J,T}\right)^{1/2}$ and each operator $\TT_i$, \ $i=1,\ldots, n$,  is uniquely defined by the relation
$$
\left( P_{\cN_J\otimes \cD_T}|_{\HH_{J,T}}\right) \TT_i^*x= (B_i^*\otimes I_{\cD_T})\left( P_{\cN_J\otimes \cD_T}|_{\HH_{J,T}}\right)x, \quad x\in \HH_{J,T},
$$
where $P_{\cN_J\otimes \cD_T}|_{\HH_{J,T}}$ is a one-to-one operator,
$ P_{\cN_J\otimes \cD_T}$ is the orthogonal projection of the Hilbert space
$\left(\cN_J\otimes \cD_T\right)\oplus \overline{\Delta_{J,T}(\cN_J\otimes \cD_{T^*})}$ onto the subspace $\cN_J\otimes \cD_T$, and  $B_1,\ldots, B_n$ are the constrained left creation operators determined by $J$.

Moreover,  $T$ is a   constrained pure row contraction if and only
if the constrained characteristic function  $\Theta_{J,T}$ is an
inner multi-analytic operator with respect to $B_1,\ldots, B_n$.
  In this case,   $T$ is unitarily equivalent to the   row
contraction
\begin{equation}\label{HH}
\left[P_
{\HH_{J,T}} (B_1\otimes I_{\cD_T})|\HH_{J,T},\ldots,
P_{\HH_{J,T}} (B_n\otimes I_{\cD_T})|\HH_{J,T}\right],
\end{equation}
where $P_{\HH_{J,T}}$ is the orthogonal projection of $\cN_J\otimes \cD_T$ onto the Hilbert space
$$\HH_{J,T}:=\left(\cN_J\otimes \cD_T\right)\ominus
\Theta_{J,T}(\cN_J\otimes \cD_{T^*}).
$$
 \end{theorem}

\begin{proof}
Consider the Hilbert space \, $ \KK_{J,T}:=\left(\cN_J\otimes
\cD_T\right)\oplus \overline{\Delta_{J,T}(\cN_J\otimes \cD_{T^*})}
$\, and  define the operator $\Phi: \cN_J\otimes \cD_{T^*}\to
\KK_{J,T}$ by setting
$$
\Phi f:=\Theta_{J,T} f\oplus \Delta_{J,T} f,\quad f\in
\cN_J\otimes \cD_{T^*}.
$$
Notice that $\Phi$ is an isometry and
\begin{equation}\label{fi}
\Phi^*(g\oplus 0)=\Theta_{J,T}^*g, \quad g\in \cN_J\otimes \cD_T.
\end{equation}
Consequently, denoting by $P_{\HH_{J,T}}$ the orthogonal projection of $\KK_{J,T}$ onto the subspace $\HH_{J,T}$, we have
\begin{equation*}
\begin{split}
\|g\|^2&= \|P_{\HH_{J,T}}(g\oplus 0)\|^2+\|\Phi \Phi^*(g\oplus 0)\|^2\\
&=\|P_{\HH_{J,T}}(g\oplus 0)\|^2+\|\Theta_{J,T}^*g\|^2
\end{split}
\end{equation*}
for any $g\in \cN_J\otimes \cD_T$. On the other hand, due to Theorem \ref{J-factor},
$$
\|K_{J,T}^* g\|^2+
\|\Theta_{J,T}^*g\|^2=\|g\|^2, \quad g\in \cN_J\otimes \cD_T.
$$
Combining the above relations, we deduce that
\begin{equation}\label{K*P}
\|K_{J,T}^* g\|=\|P_{\HH_{J,T}}(g\oplus 0)\|,  \quad g\in \cN_J\otimes \cD_T.
\end{equation}
Since $[T_1,\ldots, T_n]$ is  a constrained c.n.c. row contraction,
Corollary \ref{Co} shows that $K_{J,T}$ is  a one-to-one operator
 and $\text{\rm range}\, K_{J,T}^*$ is dense in $\cH$.

Let $x\in \HH_{J,T}$ and assume that $x\perp P_{\HH_{J,T}}(g\oplus 0)$ for any $g\in \cN_J\otimes \cD_T$.
Using the definition of $\HH_{J,T}$ and the fact that
$$
\KK_{J,T}=\left\{g\oplus 0:\ g\in \cN_J\otimes \cD_T\right\}\bigvee
\left\{\Theta_{J,T} f\oplus \Delta_{J,T} f,\quad f\in
\cN_J\otimes \cD_{T^*}\right\},
$$
we deduce that $x=0$. This shows that
$$
\HH_{J,T}=\left\{P_{\HH_{J,T}}(g\oplus 0):\ g\in \cN_J\otimes \cD_T\right\}^{-}.
$$
Hence, and due to relation \eqref{K*P},  there is a unique unitary operator
$\Gamma:\cH\to \HH_{J,T}$ such that
\begin{equation}\label{Ga}
\Gamma(K_{J,T}^* g)=P_{\HH_{J,T}}(g\oplus 0), \quad  g\in \cN_J\otimes \cD_T.
\end{equation}
Using  Theorem \ref{J-factor}, relation \eqref{fi},  and the fact that $\Phi$ is an isometry, we have
\begin{equation*}
\begin{split}
P_{\cN_J\otimes \cD_T} \Gamma K_{J,T}^* g&=
P_{\cN_J\otimes \cD_T} P_{\HH_{J,T}}(g\oplus 0)\\
&=
g-P_{\cN_J\otimes \cD_T} \Phi \Phi^*(g\oplus 0)\\
&=g-\Theta_{J,T}
\Theta_{J,T}^* g\\
&=K_{J,T} K_{J,T}^*g
\end{split}
\end{equation*}
for any $g\in \cN_J\otimes \cD_T$. Consequently, since the range of $K_{J,T}^*$ is dense in $\cH$, we deduce that
\begin{equation}
\label{PGK}
P_{\cN_J\otimes \cD_T} \Gamma=K_{J,T}.
\end{equation}

For each $i=1,\ldots, n$, let $\TT_i:\HH_{J,T}\to \HH_{J,T}$ be
defined by $\TT_i:=\Gamma T_i\Gamma^*$,\ $i=1,\ldots, n$. Due to
relation \eqref{PGK} and taking into account that the constrained
Poisson kernel $K_{J,T}$ is one-to-one, we deduce that
\begin{equation}
\label{PKGa}
P_{\cN_J\otimes \cD_T} |_{\HH_{J,T}}=K_{J,T} \Gamma^*
\end{equation}
is a one-to-one operator acting from $\HH_{J,T}$ to $\cN_J\otimes
\cD_T$. Notice also that, using relation \eqref{PKGa} and Corollary
\ref{Co}, we have
\begin{equation*}
\begin{split}
\left(P_{\cN_J\otimes \cD_T} |_{\HH_{J,T}}\right) \TT_i^*\Gamma h&=
\left(P_{\cN_J\otimes \cD_T} |_{\HH_{J,T}}\right) \Gamma T_i^* h
=K_{J,T} T_i^*h\\
&=
\left( B_i^*\otimes I_{\cD_T}\right) K_{J,T}h\\
&= \left( B_i^*\otimes I_{\cD_T}\right)
\left(P_{\cN_J\otimes \cD_T} |_{\HH_{J,T}}\right)\Gamma h
\end{split}
\end{equation*}
for any $h\in \cH$.
Hence, we deduce that

\begin{equation}
\label{def}
\left( P_{\cN_J\otimes \cD_T}|_{\HH_{J,T}}\right) \TT_i^*x= (B_i^*\otimes I_{\cD_T})\left( P_{\cN_J\otimes \cD_T}|_{\HH_{J,T}}\right)x, \qquad x\in \HH_{J,T}.
\end{equation}
Since the operator $P_{\cN_J\otimes \cD_T}|_{\HH_{J,T}}$
is one-to-one (see \eqref{PKGa}), the relation
\eqref{def}
uniquely determines the operators $\TT_i^*$, \ $i=1,\ldots, n$.

To prove that last part of the theorem,
assume that  $T:=[T_1,\ldots, T_n]$ is a   constrained pure
row contraction.
 According to Corollary \ref{Co}, the constrained Poisson kernel $K_{J,T}:\cH\to \cN_J\otimes \cD_T$
 is an isometry.
Consequently, $K_{J,T}K_{J,T}^*$ is the orthogonal projection of $\cN_J\otimes \cD_T$ onto
$K_{J,T}\cH$. According to  Theorem \ref{J-factor}, relation \eqref{J-fa} shows that $K_{J,T}K_{J,T}^*$ and
$\Theta_{J,T}\Theta_{J,T}^*$ are mutually orthogonal projections such that
$$
K_{J,T}K_{J,T}^*+\Theta_{J,T}\Theta_{J,T}^*=
I_{\cN_J\otimes \cD_T}.
$$
Therefore, $\Theta_{J,T}$ is a partial isometry, i.e., an inner multi-analytic operator and
$\Theta_{J,T}^* \Theta_{J,T}$ is a projection.
This implies that $\Delta_{J,T}$ is the projection  on
the orthogonal complement of $\text{\rm range}\, \Theta_{J,T}^*$.

Now, notice that a vector $u\oplus v\in \KK_{J,T}$ is in $\HH_{J,T}$ if and only if
$$
\left<
u\oplus v,\Theta_{J,T}f\oplus \Delta_{J,T}f\right>=0\quad
\text{ for any } \ f\in \cN_J\otimes \cD_{T^*}.
$$
This is equivalent to
\begin{equation}
\label{TH*}
\Theta_{J,T}^*u+\Delta_{J,T}v=0.
\end{equation}
Due to the above observations, we have $\Theta_{J,T}^*u\perp \Delta_{J,T}v$. This shows that
relation \eqref{TH*} holds if and only if
$\Theta_{J,T}^*u=0$ and $v=0$. Consequently,
$$\HH_{J,T}=\left(\cN_J\otimes \cD_T\right)\ominus
\Theta_{J,T}(\cN_J\otimes \cD_{T^*}).
$$
In this case, $P_{\cN_J\otimes \cD_T}|_{\HH_{J,T}}$ is the restriction operator and relation \eqref{def} implies
$$
\TT_i=
P_
{\HH_{J,T}} (B_i\otimes I_{\cD_T})|\HH_{J,T}, \quad i=1,\ldots, n.
$$
Conversely,  if $\Theta_{J,T}$ is inner, then it is a partial
isometry. Theorem \ref{J-factor} implies that  $K_{J,T}$ is a
partial isometry. Since $T$ is c.n.c., Corollary \ref{Co} implies $
\text{\rm SOT-}\,\lim_{k\to\infty}\sum_{|\alpha|=k} T_\alpha
T_\alpha^*=0, $ which proves that $T$ is a pure row contraction.
    This completes the proof.
\end{proof}

As in the noncommutative case \cite{Po-charact}, one can easily prove the following.
\begin{proposition}
If $T:=[T_1,\ldots, T_n]$, $T_i\in B(\cH)$, is a  $J$-constrained c.n.c. row contraction, then $\Theta_{J,T}$ is outer if and only if there is no element $h\in \cH$, $h\neq 0$, such that
 $
 \lim_{k\to\infty}\sum_{|\alpha|=k} T_\alpha T_\alpha^*h=0.
 $
\end{proposition}
\begin{proof}
Due to Corollary \ref{Co}, the condition  above   is equivalent to
$\ker \left( I-K_{J.T}^* K_{J,T}\right)=\{0\}$. Using Theorem
\ref{J-factor}, we deduce that the latter equality is equivalent to
$$
\ker \Theta_{J,T} \Theta_{J,T}^*=\ker \left( I-K_{J.T} K_{J,T}^*\right)=\{0\},
$$
which is the same as  $\Theta_{J,T}$ having dense range.
The proof is complete.
\end{proof}

\begin{remark} If $J=\{0\}$ in Theorem $\ref{model}$, one can recover the model theorem for arbitrary c.n.c. row contractions \cite{Po-charact}.
\end{remark}

Let  $\Phi\in \cW(W_1,\ldots, W_n)\bar\otimes B(\cK_1, \cK_2)$
and
$\Phi'\in \cW(W_1,\ldots, W_n)\bar\otimes B(\cK_1', \cK_2')$ be two multi-analytic operators with respect to $B_1,\ldots, B_n$. We say that $\Phi$ and $\Phi'$ coincide
if there are two unitary operators $\tau_j\in B(\cK_j, \cK_j')$, $j=1,2$,  such that
$$
\Phi'(I_{\cN_J}\otimes \tau_1)=(I_{\cN_J}\otimes \tau_2) \Phi.
$$
We remark that if $1\in \cN_J$, then the $C^*$-algebra
$C^*(B_1,\ldots, B_n)$ is irreducible (see \cite{Po-constrained}). In this case,   the operators $\Phi$ and $\Phi'$ coincide
if and only if  there are two multi-analytic operators
$U_j: \cN_J\otimes \cK_j\to \cN_J\otimes\cK_j'$
such that   $\Phi' U_1=U_2\Phi$.

The next result shows that the constrained characteristic function is a complete unitary invariant for c.n.c. constrained row contractions.

\begin{theorem}\label{u-inv}
Let $J\neq F_n^\infty$ be a WOT-closed two-sided ideal of
$F_n^\infty$ and let $T:=[T_1,\ldots, T_n]$, \ $T_i\in B(\cH)$, and
$T':=[T_1',\ldots, T_n']$,\ $T_i'\in B(\cH')$, be two
$J$-constrained  c.n.c. row contractions. Then $T$ and $T'$ are
unitarily equivalent if and only if their  constrained
characteristic functions $\Theta_{J,T}$  and $\Theta_{J,T'}$
coincide.
\end{theorem}
\begin{proof}
Assume that $T$ and $T'$ are unitarily equivalent and let $U:\cH\to \cH'$ be a unitary operator such that
$T_i=U^*T_i'U$ for any $i=1,\ldots, n$. Simple computations reveal that
$$
U\Delta_T=\Delta_{T'}U \quad \text{ and }\quad
(\oplus_{i=1}^n U)\Delta_{T^*}=\Delta_{T'^*}(\oplus_{i=1}^n U).
$$
Define the unitary operators $\tau$ and $\tau'$ by
setting
$$\tau:=U|\cD_T:\cD_T\to \cD_{T'} \quad \text{ and }\quad
\tau':=(\oplus_{i=1}^n U)|\cD_{T^*}:\cD_{T*}\to \cD_{T'^*}.
$$
Taking into account the definition of the constrained characteristic function, it is easy to see that
$$
(I_{\cN_J}\otimes \tau)\Theta_{J,T}=\Theta_{J,T'}(I_{\cN_J}\otimes \tau').
$$

Conversely, assume that the constrained characteristic functions  of $T$ and $T'$ coincide.  Then  there exist unitary operators
$\tau:\cD_T\to \cD_{T'}$ and $\tau_*:\cD_{T^*}\to \cD_{{T'}^*}$ such that
\begin{equation}\label{com}
(I_{\cN_J}\otimes \tau)\Theta_{J,T}=\Theta_{J,T'}(I_{\cN_J}\otimes \tau_*).
\end{equation}
It is easy to see that \eqref{com} implies \, $
\Delta_{J,T}=\left(I_{\cN_J}\otimes \tau_*\right)^*
\Delta_{J,T'}\left(I_{\cN_J}\otimes \tau_*\right) $ and
$$
\left(I_{\cN_J}\otimes \tau_*\right)\overline{\Delta_{J,T}(\cN_J\otimes \cD_{T^*})}=
\overline{\Delta_{J,T'}(\cN_J\otimes \cD_{{T'}^*})}.
$$
 Using the notations of Theorem \ref{model}, we define the unitary operator $U:\KK_{J,T}\to \KK_{J,T'}$
by setting $U:=(I_{\cN_J}\otimes \tau)\oplus (I_{\cN_J}\otimes
\tau_*)$.  Straightforward computations reveal that the operator
$\Phi:\cN_J\otimes \cD_{T^*}\to \KK_{J,T}$, defined in the proof of
Theorem \ref{model}, and the corresponding $\Phi'$ satisfy the
relations
\begin{equation}
\label{Uni1}
U \Phi\left(I_{\cN_J}\otimes \tau_*\right)^*=\Phi'
\end{equation}
and
\begin{equation}
\label{Uni2}
\left(I_{\cN_J}\otimes \tau\right) P_{\cN_J\otimes \cD_T}^{\KK_{J,T}} U^*=P_{\cN_J\otimes \cD_{T'}}^{\KK_{J,T'}},
\end{equation}
where $P_{\cN_J\otimes \cD_T}^{\KK_{J,T}}$ is the orthogonal projection of $\KK_{J,T}$ onto $\cN_J\otimes \cD_T$.
Notice that relation \eqref{Uni1} implies
\begin{equation*}
\begin{split}
U\HH_{J,T}&=U\KK_{J,T}\ominus U\Phi(\cN_J\otimes \cD_{T^*})\\
&=\KK_{J,T'}\ominus \Phi'(I_{\cN_J}\otimes \tau_*)(\cN_J\otimes \cD_{T^*})\\
&=\KK_{J,T'}\ominus \Phi' (\cN_J\otimes \cD_{{T'}^*}).
\end{split}
\end{equation*}
Therefore, the operator $U|_{\HH_{J,T}}:\HH_{J,T}\to \HH_{J, T'}$ is
unitary. On the other hand, we have
\begin{equation}
\label{intertw}
(B_i^*\otimes I_{\cD_{T'}})(I_{\cN_J}\otimes \tau)=
(I_{\cN_J}\otimes \tau)(B_i^*\otimes I_{\cD_{T}}).
\end{equation}
Now, let $\TT:=[\TT_1,\ldots \TT_n]$ and $\TT':=[\TT_1',\ldots \TT_n']$ be the models provided by Theorem \ref{model}  for the row contractions $T$ and $T'$, respectively.
Using the relation \eqref{def}  for $T'$ and $T$,  as well as relations  \eqref{Uni2} and \eqref{intertw}, we deduce  that
\begin{equation*}
\begin{split}
P_{\cN_J\otimes \cD_{T'}}^{\KK_{J,T'}}{\TT_i'}^*Ux &=  (B_i^*\otimes
I_{\cD_{T'}}) P_{\cN_J\otimes \cD_T}^{\KK_{J,T}}Ux
=(B_i^*\otimes I_{\cD_{T'}})(I_{\cN_J}\otimes \tau) P_{\cN_J\otimes \cD_T}^{\KK_{J,T}}x\\
&=(I_{\cN_J}\otimes \tau)(B_i^*\otimes I_{\cD_{T}})P_{\cN_J\otimes
\cD_T}^{\KK_{J,T}}x =(I_{\cN_J}\otimes \tau) P_{\cN_J\otimes
\cD_T}^{\KK_{J,T}}
\TT_i^*x\\
&= P_{\cN_J\otimes \cD_{T'}}^{\KK_{J,T'}}U \TT_i^*x
\end{split}
\end{equation*}
for any $x\in \HH_{J,T}$ and $i=1,\ldots, n$. Consequently, since
$P_{\cN_J\otimes \cD_{T'}}^{\KK_{J,T'}}|_{\HH_{J,T'}}$ is a
one-to-one operator (see  relation \eqref{PKGa}), we obtain
$$
\left(U|_{\HH_{J,T}}\right) \TT_i^*={\TT_i'}^*\left(U|_{\HH_{J,T}}\right),\quad i=1,\ldots,n.
$$
 Now, using Theorem \ref{model}, we conclude that $T$ and $T'$ are unitarily equivalent. The proof is complete.
\end{proof}

As in the noncommutative case \cite{Po-curvature}, one can prove the following.
\begin{proposition}
Let $J\neq F_n^\infty$ be a WOT-closed two -sided ideal of $F_n^\infty$ such that $1\in \cN_J$. If $T:=[T_1,\ldots, T_n]$ is a  $J$-constrained c.n.c.  row contraction, then $T$ is unitarily equivalent to
a  constrained shift $[B_1\otimes I_\cK,\ldots, B_n\otimes I_\cK]$ for some Hilbert space $\cK$ if and only if $\Theta_{J,T}=0$.
\end{proposition}
\begin{proof}
If $T=[B_1\otimes I_\cK,\ldots, B_n\otimes I_\cK]$, then
$K_{J,T}f=f$ for $f\in \cN_J\otimes \cK$. Indeed, since $1\in \cN_J$,  for any
$f=\sum\limits_{\alpha\in \FF_n^+} e_\alpha\otimes k_\alpha$  in $N_J\otimes \cK\subseteq F^2(H_n)\otimes \cK$, we have

\begin{equation*}
\begin{split}
K_{J,T}f&=\sum\limits_{\alpha\in \FF_n^+} P_{\cN_J}e_\alpha\otimes P_\cK(B_\alpha^*\otimes I_\cK)f=
 \sum\limits_{\alpha\in \FF_n^+} P_{\cN_J}e_\alpha\otimes P_\cK(S_\alpha^*\otimes I_\cK)f\\
&=
\sum\limits_{\alpha\in \FF_n^+} P_{\cN_J}e_\alpha\otimes k_\alpha=P_{\cN_J\otimes \cK} f=f.
 \end{split}
 \end{equation*}
Now, Theorem \ref{J-factor} shows that $\Theta_{J,T}=0$.
Conversely, if $\Theta_{J,T}=0$, then Theorem \ref{model} shows that $T$ is unitarily
equivalent to
the   constrained shift $[B_1\otimes I_{\cD_T},\ldots, B_n\otimes I_{\cD_T}]$.
\end{proof}

\begin{remark} All the  results of this section can be written  in the particular
case when $T:=[T_1,\ldots,T_n]$ is a
$q$-commuting c.n.c. row contractions, i.e.,
$$
T_iT_j=q_{ji}T_jT_i,\quad 1\leq i<j\leq n,
$$
where $q_{ij}\in \CC$.
 \end{remark}
Notice that  $T$ is  $q$-commuting
if and only if it is a $J_q$-constrained row contraction,
where
$J_q$ is  the WOT-closed two-sided ideal of $F_n^\infty$ generated by
the $q$-commutators $S_iS_j-q_{ji}S_jS_i$, $1\leq i<j\leq n$, where
 $S_1,\ldots, S_n$ are the left creation operators  on the full Fock space.
We refer to \cite{BB} for related results on $q$-commuting row
contractions.

 We remark that in the particular case when $q_{ij}=1$
we obtain a model theory for commuting c.n.c. row contractions.

       %

      \end{document}